\documentclass [11pt,oneside]{article}

\usepackage{amssymb}
\usepackage{amsfonts}
\usepackage{amsmath}
\usepackage{amsthm}
\usepackage{epsfig}

\newtheorem{lemma}{Lemma}[section]
\newtheorem{teo}[lemma]{Theorem}
\newtheorem{rem}[lemma]{Remark}
\newtheorem{prop}[lemma]{Proposition}
\newtheorem{cor}[lemma]{Corollary}

\newcommand{\ordtred}{\matD^{\, 3}_{\rm o}}
\newcommand{\cyctred}{\matD^{\, 3}_{\rm c}}
\newcommand{\vertred}{\matD^{\, 3}_{\rm v}}
\newcommand{\orddues}{\matS^2_{\rm o}}
\newcommand{\cycdues}{\matS^2_{\rm c}}
\newcommand{\verdues}{\matS^2_{\rm v}}
\newcommand{\ordtres}{\matS^3_{\rm o}}
\newcommand{\cyctres}{\matS^3_{\rm c}}
\newcommand{\vertres}{\matS^3_{\rm v}}

\newcommand{\matN} {\ensuremath {\mathbb{N}}}
\newcommand{\matR} {\ensuremath {\mathbb{R}}}

\newcommand{\matC} {\ensuremath {\mathbb{C}}}
\newcommand{\matP} {\ensuremath {\mathbb{P}}}

\newcommand{\matD} {\ensuremath {\mathbb{D}}}
\newcommand{\matS} {\ensuremath {\mathbb{S}}}

\newcommand{\calX} {\ensuremath {\mathcal{X}}}

\newcommand{\calS} {\ensuremath {\mathcal{S}}}

\newcommand{\nota} [1] {\caption{\footnotesize{#1}}}

\newfont{\Got}{eufm10 scaled 1200}

\font\titsc=cmcsc10 scaled 1200

\newcommand{\dimo}[1]{\vspace{2pt}\noindent\textit{Proof of \ref{#1}}.\ }
\newcommand{\finedimo}{{\hfill\hbox{$\square$}\vspace{2pt}}}
\newcommand{\mettifig}[1]{\epsfig{file=#1}}

\author{Carlo \titsc{Petronio}\thanks{Research supported
by the INTAS Project ``CalcoMet-GT'' 03-51-3663}}

\title{Complexity of 3-orbifolds}

\begin{document}

\maketitle

\begin{abstract}
    \noindent We extend Matveev's theory of complexity for 3-manifolds,
    based on simple spines, to (closed, orientable, locally orientable) 3-orbifolds.
    We prove naturality and finiteness for irreducible 3-orbifolds, and, with certain
    restrictions and subtleties, additivity under orbifold connected sum.
    We also develop the theory
    of handle decompositions for $3$-orbifolds and the corresponding theory of
    normal $2$-suborbifolds.
  \vspace{4pt}

\noindent MSC (2000): 57M99.
\end{abstract}

\section*{Introduction}
The aim of this paper is to extend Matveev's theory of complexity from
$3$-manifolds to $3$-orbifolds, and in particular, after giving the appropriate definition, to prove
the following properties, established for manifolds in~\cite{Matveev:AAM}:
\begin{itemize}
\item Naturality and finiteness for irreducible orbifolds;
\item Additivity under orbifold connected sum.
\end{itemize}
To quickly summarize our results, we anticipate that we have fully
achieved the former of these tasks, with the most natural
definition of complexity. Concerning the latter task, we have
shown additivity to hold in a strict (and subtle) fashion for two
of the three possible types of connected sum, and we have
established a two-sided linear estimate for the third type.

\vspace{.3cm}

We now give a more detailed account of our results. We always
consider closed, orientable, locally orientable $3$-orbifolds. To
define the \emph{complexity} $c(X)$ of such an $X$ we consider the
set of simple polyhedra $P$ embedded in $|X|$ so that $P$ is
transversal to the singular set $S(X)$ and $X\setminus P$ consists
of disjoint open discal $3$-orbifolds. Then we take the minimum
over all $P$ of the number of vertices of $P$ plus a contribution
for each intersection point between $P$ and $S(X)$, the
contribution being $p-1$ if the order of $S(X)$ at the point is
$p$. We can then state the naturality and finiteness properties
proved below:

\begin{teo}\label{naturality:intro:teo}
With certain well-understood exceptions, if $P$ realizes the complexity of
an irreducible orbifold $X$ then $P$ is a special polyhedron, and dual to $P$ there is
a triangulation of $X$. The number of exceptional orbifolds of each given complexity is finite.
\end{teo}

\begin{cor}\label{finiteness:intro:cor}
For each positive integer $n$ the set of
irreducible orbifolds $X$ such that $c(X)\leqslant n$
is finite.
\end{cor}

Turning to additivity, first recall that the operation of connected sum of two manifolds or orbifolds
consists in making a puncture in each of the given objects and then gluing the
resulting boundary spheres. For manifolds, there is only one way to make the puncture, so the
operation of connected sum is uniquely defined (at least in a connected and oriented context).
In the analogue operation for orbifolds,
one has to distinguish according to the nature of the point at which
the puncture is made, which can be non-singular, singular but not a vertex of the singular
locus, or a vertex of the singular locus. Depending on this nature, the operation of
connected sum is called of \emph{ordinary}, \emph{cyclic}, or \emph{vertex} type.
Each type of operation has a corresponding identity element, given by the so-called
\emph{ordinary}, \emph{cyclic}, and \emph{vertex spherical $3$-orbifold}.
A connected sum of a certain type with the spherical $3$-orbifold of the same type will
be called \emph{trivial}. The following result, proved in~\cite{orb:split},
is necessary to state our results on additivity.

\begin{teo}\label{splitting:cite:teo}
Let $X$ be a closed locally orientable $3$-orbifold. Suppose that $X$ does not contain
any bad $2$-suborbifold, and that every spherical $2$-suborbifold of $X$ is separating.
Then $X$ can be realized as connected sum of irreducible $3$-orbifolds in such a way that,
even after reordering the sums, there is no trivial sum.
Any two such realizations involve the same irreducible summands and the same types of sums
(including orders). Moreover, if the realization of $X$ does not involve vertex connected sums,
for all $p\geqslant 2$ the number $\nu(p)$ of
$p$-cyclic connected sums involving at least one singular
component without vertices is independent of the realization, so it is a function of $X$ only.
\end{teo}

A realization of $X$ as described in the previous statement will
be called an \emph{efficient connected sum}. The condition that
there should not be trivial sums \emph{even after reordering} is
due to the existence of phenomena of the following type. If we
first take the ordinary connected sum of some non-singular $X$
with a cyclic spherical $3$-orbifold, and then we take the cyclic
connected sum of the result with some other $Y$, then both sums,
taken in this order, are non-trivial, but the result is just the
ordinary connected sum of $X$ and $Y$, so the two-step sum is
obviously inefficient. We now have our main statement about
additivity:

\begin{teo}\label{additivity:intro:teo}
Let a $3$-orbifold $X$ as in Theorem~\ref{splitting:cite:teo}
be the efficient connected sum of irreducible orbifolds $X_1,\ldots,X_n$
without vertex connected sums. For $p\geqslant 2$ let $\nu(p)$ be
as in Theorem~\ref{splitting:cite:teo}. Then
$$c(X)=\sum_{j=1}^{n-1} c(X_i)-\sum_{p=2}^\infty \nu(p)\cdot (p-1).$$
\end{teo}

When the realization of $X$ includes vertex connected sums we only have two-sided
linear estimates on the complexity of $X$ in terms of the complexity of its summands:

\begin{teo}\label{estimates:intro:teo}
Let a $3$-orbifold $X$ as in Theorem~\ref{splitting:cite:teo}
be the efficient connected sum of irreducible orbifolds $X_1,\ldots,X_n$. Then
$$\frac 1{4^{n-1}}\cdot \big(c(X_1)+\ldots+c(X_n)\big)\leqslant
c(X)\leqslant 6^{n-1}\cdot \big(c(X_1)+\ldots+c(X_n)\big).$$
\end{teo}

We mention that to establish the additivity properties of complexity we have
developed an apparently new portion of the general theory of orbifolds, namely the notion of
handle decomposition. The reader will find the details on this notion in
Section~\ref{handles:section}, together with the analogue in the orbifold
setting of Haken's theory of normal surfaces with respect to handle decompositions.

As a motivation for our interest in $3$-orbifolds, we would like
to recall here that orbifolds play a central r\^ole in
$3$-dimensional geometric topology, and in particular in
Thurston's geometrization program, see~\cite{BLP}.

\section{Preliminaries and main definitions}\label{def:section}
In this section we briefly recall the notions of orbifold and of simple polyhedron,
and we define the complexity of a $3$-orbifold.

\paragraph{Local structure of orbifolds}
We will not cover here the general theory of orbifolds, referring the reader to the
milestone~\cite{thurston:notes}, to the excellent and very recent~\cite{BMP},
and to the comprehensive bibliography of the latter.
We just recall that an orbifold of dimension $n$ is a topological space with a singular smooth
structure, locally modelled on a quotient of $\matR^n$ under the action of a finite
group of diffeomorphisms. We will only need to refer to the cases $n=2$ and $n=3$, and we
will confine ourselves to orientation-preserving diffeomorphisms. In addition, all our
orbifolds will be compact and globally orientable.

Under these assumptions one can see that a 2-orbifold $\Sigma$
consists of a compact orientable \emph{support} surface $|\Sigma|$ together with a finite
collection $S(\Sigma)$ of points in the interior of $|\Sigma|$, the \emph{cone} points,
each carrying a certain \emph{order} in $\{p\in\matN:\ p\geqslant2\}$.

Analogously, a $3$-orbifold $X$ is given by a compact \emph{support} 3-manifold $|X|$ together
with a \emph{singular set} $S(X)$. Here $S(X)$ is a finite collection of circles and unitrivalent
graphs tamely embedded in $|X|$, with the univalent vertices
given by the intersection with $\partial|X|$. Moreover each component of $S(X)$ minus the vertices
carries an \emph{order} in $\{p\in\matN:\ p\geqslant2\}$, with the restriction
that the three germs of edges incident to each vertex should
have orders $(2,2,p)$, for arbitrary $p$, or $(2,3,p)$, for $p\in\{3,4,5\}$.

\paragraph{Bad, spherical, and discal orbifolds}
An \emph{orbifold-covering} is a map between orbifolds locally modelled on a map of
the form $\matR^n/_\Delta\to\matR^n/_\Gamma$, naturally defined whenever
$\Delta<\Gamma<{\rm Diff}_+(\matR^n)$. An orbifold is called \emph{good} when
it is orbifold-covered by a manifold, and \emph{bad} when it is not good.
In the sequel we will need the following easy result:

\begin{lemma}\label{bad:lem}
The only bad closed $2$-orbifolds are $(S^2;p)$, the $2$-sphere with one cone point of
order $p$, and $(S^2;p,q)$, the $2$-sphere with cone points of orders $p\neq q$.
\end{lemma}

We now introduce some notation and terminology repeatedly
used below. We define $\ordtred$ to be $D^3$, the \emph{ordinary discal} $3$-orbifold,
$\cyctred(p)$ to be $D^3$ with singular set
a trivially embedded arc with arbitrary order $p$,
and $\vertred(p,q,r)$ to be $D^3$ with singular set
a trivially embedded
``Y-graph'' with edges of admissible orders $p,q,r$.
We will call $\cyctred(p)$ and $\vertred(p,q,r)$ respectively
\emph{cyclic discal} and \emph{vertex discal} $3$-orbifolds, and we will employ
the shortened notation $\cyctred$ and $\vertred$ to denote
cyclic and vertex discal $3$-orbifolds with generic orders.

We also define the \emph{ordinary, cyclic}, and \emph{vertex spherical} $2$-orbifolds,
denoted respectively by $\orddues$, $\cycdues(p)$, and $\verdues(p,q,r)$,
as the 2-orbifolds bounding the corresponding discal 3-orbifolds
$\ordtred$, $\cyctred(p)$, and $\vertred(p,q,r)$.
Finally, we introduce the \emph{ordinary, cyclic}, and \emph{vertex spherical} $3$-orbifolds,
denoted respectively by $\ordtres$, $\cyctres(p)$, and $\vertres(p,q,r)$,
as the 3-orbifolds obtained by mirroring the corresponding discal 3-orbifolds
$\ordtred$, $\cyctred(p)$, and $\vertred(p,q,r)$ in their boundary.
The spherical 2- and 3-orbifolds with generic orders will be denoted by $\matS^2_*$ and $\matS^3_*$.

\paragraph{2-suborbifolds and irreducible 3-orbifolds}
We say that a $2$-orbifold $\Sigma$ is a \emph{suborbifold} of a $3$-orbifold $X$ if $|\Sigma|$
is embedded in $|X|$ so that $|\Sigma|$ meets $S(X)$ transversely (in particular,
it does not meet the vertices), and $S(\Sigma)$ is given precisely by $|\Sigma|\cap S(X)$,
with matching orders.

A spherical $2$-suborbifold $\Sigma$ of a $3$-orbifold $X$ is
called \emph{essential} if it does not bound in $X$ a discal
$3$-orbifold. A $3$-orbifold $X$ is called \emph{irreducible} if
it does not contain any bad 2-suborbifold and every spherical
2-suborbifold of $X$ is inessential (in particular, it is
separating).

\begin{figure}
    \begin{center}
    \mettifig{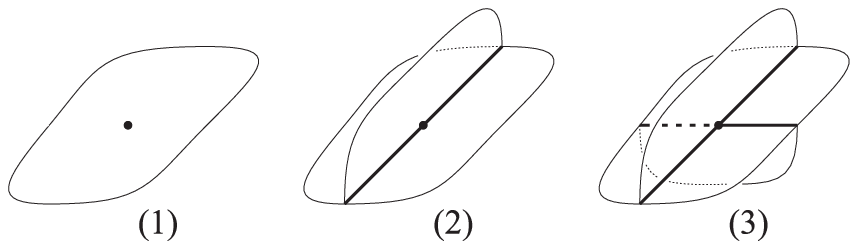,width=8cm}
    \nota{Local aspect of an almost-special polyhedron.} \label{almostspecial:fig}
    \end{center}
    \end{figure}

\paragraph{Simple and special polyhedra}
From now on we will employ the piecewise linear viewpoint, which
is equivalent to the smooth one in dimensions $2$ and $3$. We will
use the customary notions of PL topology without recalling them,
see~\cite{RS}. A \emph{simple polyhedron} is a compact polyhedron
$P$ such that the link of each point of $P$ can be embedded in the
space given by a circle with three radii. In particular, $P$ has
dimension at most $2$. Finite graphs and closed surfaces are
examples of simple polyhedra. A point of a simple polyhedron is
called a \emph{vertex} if its link is precisely given by a circle
with three radii. A regular neighbourhood of a vertex is shown in
Fig.~\ref{almostspecial:fig}-(3).
From the figure one sees that the vertices are isolated, whence
finite in number.  Graphs and surfaces do not contain vertices.
The complexity $c(P)$ of $P$ is the number of vertices that $P$ contains.

Two more restrictive types of polyhedra will be used below.
A simple polyhedron $P$ is called \emph{almost-special} if
the link of each point of $P$ is given by a circle with either zero, or two, or three radii.
The local aspects of $P$ are correspondingly shown in Fig.~\ref{almostspecial:fig}.
The points of type (2) or (3) are called \emph{singular}, and the set
of singular points of $P$ is denoted by $S(P)$. We will say that $P$ is \emph{special}
if it is almost-special, $S(P)$ contains no circle component, and $P\setminus S(P)$
consists of open $2$-discs.

\paragraph{Duality}
The following duality result, that we will state in the closed
context only, is well-known. We call \emph{triangulation} of a
3-manifold $M$ a realization of $M$ as a gluing of a finite number
of tetrahedra along a complete system of simplicial pairings of
the lateral faces. Note that we allow multiple and
self-adjacencies of the tetrahedra, thus relaxing the traditional
requirements for a triangulation used in PL topology.

\begin{figure}
    \begin{center}
    \mettifig{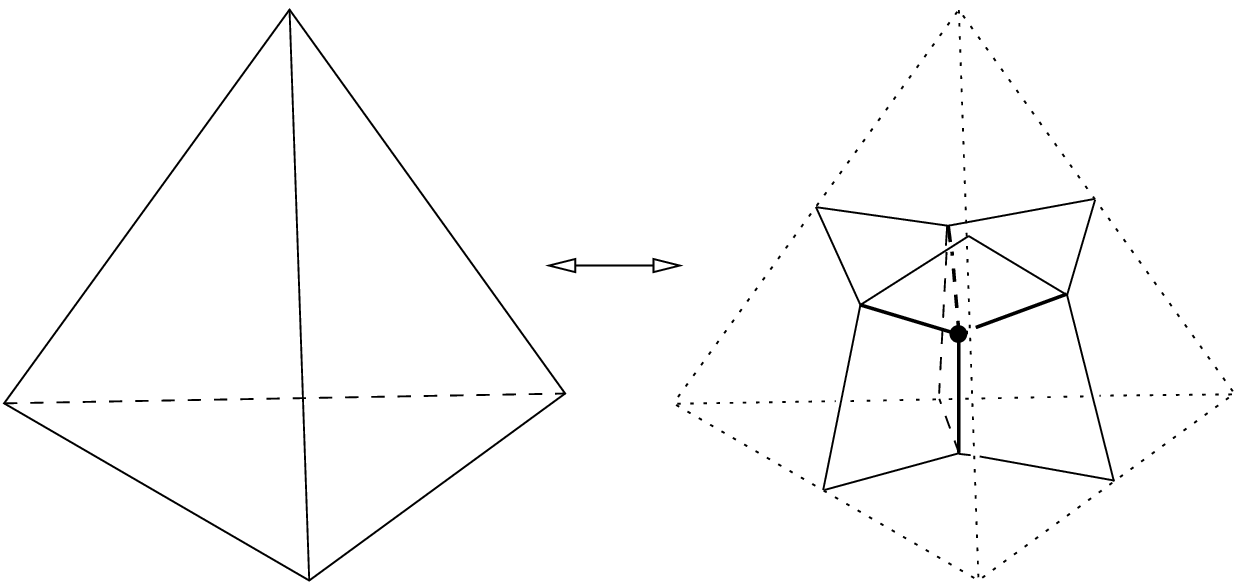,width=8cm}
    \nota{Duality between triangulations and special polyhedra.} \label{duality:fig}
    \end{center}
    \end{figure}

\begin{prop}\label{duality:basic:prop}
The set of triangulations of a $3$-manifold $M$ corresponds
bijectively to the set of special polyhedra $P$ embedded in $M$ in
such a way that $M\setminus P$ is a union of open $3$-discs. Both
triangulations and special polyhedra are viewed up to isotopy, and
the polyhedron corresponding to a triangulation is the
$2$-skeleton of the dual cellularization, as shown in
Fig.~\ref{duality:fig}.
\end{prop}

\paragraph{Spines and complexity}
Let $X$ be a (closed, orientable, locally orientable) 3-orbifold. Let $P$ be a simple polyhedron
contained in $|X|$. We will say that $P$ is a \emph{spine} of $X$ if the following holds:
\begin{itemize}
\item The intersections between $P$ and $S(X)$ occur only at
surface points of $P$ and non-vertex points of $S(X)$, and they
are transverse; \item If $U(P)$ is an open regular neighbourhood
of $P$ in $|X|$, each component of $X\setminus U(P)$ is isomorphic
to one of the discal $3$-orbifolds $\matD^{\, 3}_*$.
\end{itemize}
In the case of a manifold, \emph{i.e.} for $S(X)=\emptyset$, we
recover Matveev's condition that $X\setminus P$ consists of open
$3$-discs, \emph{i.e.} that $X$ minus some points collapses onto
$P$. Note that, as opposed to the manifold case, the existence of
a non-empty singular locus in an orbifold forces a spine to have
some 2-dimensional stratum. From now on we will identify
$X\setminus U(P)$ to $X\setminus P$, making a distinction only
when it really matters. The reader will easily check that this
choice does not cause any ambiguity.

If $P$ is a spine of $X$ as above, we define the following function,
which depends not only on $P$ as an abstract polyhedron, but also on
its embedding in $|X|$ relative to $S(X)$:
$$c(P,S(X))=c(P)+\sum\Big\{p-1:\ x\in P\cap S(X),
\textrm{\ the\ order\ of\ }S(X)\textrm{\ at\ }x\textrm{\ is\ }p\Big\}.$$
We now define $c(X)$, the \emph{complexity} of $X$, as the minimum of $c(P,S(X))$
over all spines $P$ of $X$. We note that $c$ is always well-defined because
every orbifold has simple spines: see for instance Section~\ref{irreducible:section}, where
the relation is discussed between the special spines of a $3$-orbifold $X$
and the (suitably defined) triangulations of $X$.

We will say that a spine $P$ of $X$ is \emph{minimal} if
$c(P,S(X))=c(X)$ and every proper subpolyhedron of $P$ which is
also a spine of $X$ is actually homeomorphic to $P$.

\begin{rem}
\emph{Each of the following may appear to be a more natural condition
to include in the definition that $P$ is \emph{minimal}, besides the requirement
that $c(P,S(X))=c(X)$:
\begin{itemize}
\item[(A)] $P$ does not collapse onto any proper subpolyhedron which is also a spine of $X$;
\item[(B)] $P$ does not contain any proper subpolyhedron which is also a spine of $X$.
\end{itemize}
Condition (A) is however not satisfactory, because a minimal spine
would remain minimal after a ``bubble move'' (the gluing of a
$2$-disc along the boundary of a $2$-disc contained in $P$ and
disjoint from $S(X)$ and $S(P)$). This phenomenon already appears
in the manifold case, where Matveev uses condition (B). However
(B) indeed is not the right condition for orbifolds, for otherwise
some orbifolds, \emph{e.g.} $\cyctres(p)$, would not have minimal
spines --- see the proof of Lemma~\ref{Up:lem}. One drawback of
our definition of minimality is that it does not immediately imply
that minimal spines exist, and it will require some efforts to
prove that this is actually the case.}
\end{rem}

\section{Minimal spines of irreducible orbifolds}\label{irreducible:section}
In this section we prove that, with some well-understood
exceptions, an irreducible orbifold has special minimal spines.
This will imply that for any given $n$ there is only a finite
number of irreducible 3-orbifolds having complexity $n$, a result
which opens the way to computer enumeration of orbifolds in order
of increasing complexity, analogous to that carried out for
manifolds by Martelli, Matveev, and the author,
see~\cite{Bruno:survey}.

We also give a natural notion of triangulation for an orbifold,
and we show that a special spine is typically dual to a
triangulation, which implies that our notion of complexity is a
very natural one.

\paragraph{Singular arcs meet the spine once}
We begin with an easy result
that we will use repeatedly both in this section and in Section~\ref{additivity:section}.

\begin{lemma}\label{arc:lem}
Let $P$ be a spine of an orbifold $X$
such that $c(P,S(X))=c(X)$. Let $\alpha$ be a connected component
of $S(X)$ minus the vertices. Then $\alpha\cap P$
consists of precisely one point.
\end{lemma}

\begin{proof}
At least one intersection point exists otherwise some component of
$X\setminus P$ would contain either two vertices or a non-simply
connected portion of $S(X)$. Suppose there are there are at least
two intersection points, and let $\beta$ be the open arc of
$\alpha$ between two consecutive ones $x_0$ and $x_1$. Now $\beta$
is contained in one of the open $3$-discs $B$ of which $X\setminus
P$ consists. Let $B'$ be the component of $X\setminus P$ which
contains the portion of $\alpha$ past $x_1$. A priori we could
have $B'=B$, but $S(B)$ consists of $\beta$ only, which implies
that $B'\neq B$. Therefore we can puncture $P$ near $x_1$, getting
a new spine $P'$ of $X$. Moreover $c(P',S(X))<c(P,S(X))$, a
contradiction.
\end{proof}

\paragraph{Nuclear and minimal spines}
A simple polyhedron is called \emph{nuclear} if it cannot be
collapsed onto a proper subpolyhedron. Of course a minimal spine
of a manifold is always nuclear.  We will now show that orbifolds
admit minimal spines, and that these spines are most often
nuclear.

\begin{lemma}\label{Up:lem}
Let $X$ be any orbifold and $Y$ be the ordinary connected sum of $X$ and $\cyctres(p)$.
Then $c(Y)=c(X)+(p-1)$ , and $Y$ admits minimal spines if $X$ does.
\end{lemma}

Before proving the lemma we note that, if we apply it
to $X=\ordtres$, we deduce that $c(\cyctres(p))=p-1$.

\begin{proof}
Suppose $P$ is a spine of $Y$ with $c(P,S(Y))=c(Y)$.
Let us denote by $U_p$ the singular set of $\cyctres(p)$, an order-$p$ unknot in the $3$-sphere.
Since $|Y|=|X|$ and $S(Y)=S(X)\sqcup U_p$
then $P$ is also a spine of $X$, and
$$c(P,S(Y))\geqslant c(P,S(X))+(p-1),$$
because $U_p$ must meet $P$ in at least one point. This proves
that $c(Y)\geqslant c(X)+(p-1)$. We will now show that,
conversely, if $Q$ is a spine of $X$, we can construct a spine $P$
of $Y$ such that $c(P,S(Y))=c(Q,S(X))+(p-1)$, which will imply the
first assertion. We construct $P$ so that it will be automatically
minimal if $Q$ is minimal, which will prove the second assertion.
To do so, we distinguish according to whether there exists a
component of $X\setminus P$ which is an ordinary 3-disc. If there
is one, we construct $P$ by attaching to $Q$ a lollipop contained
in this $3$-disc, see Fig.~\ref{unknot:fig}-centre. Otherwise, we
attach to $Q$ a lollipop wrapped in a sphere, as in
Fig.~\ref{unknot:fig}-right.
    \begin{figure}
    \begin{center}
    \input{unknot.pstex_t} 
    \nota{Spine of an orbifold after insertion of an extra unknotted singular component.} \label{unknot:fig}
    \end{center}
    \end{figure}
Of course $c(P,S(Y))$ has the required value and it is easy to
show that, with the only exception where $Q$ is a point, so
$X=\ordtres$ and $Y=\cyctres(p)$, the spine $P$ constructed is
minimal if $Q$ was. A $2$-disc is of course a minimal spine of
$\cyctres(p)$, and the proof is complete.
\end{proof}

\begin{prop}\label{nuclear:prop}
Let $X$ be a $3$-orbifold such that:
\begin{itemize}
\item $X$ does not contain any bad $2$-suborbifold;
\item $X$ cannot be expressed as an ordinary connected sum between some orbifold and some $\cyctres(p)$.
\end{itemize}
Then $X$ admits minimal spines, and any such spine is nuclear.
\end{prop}

\begin{proof}
We begin by proving that if we collapse a spine $P$ of $X$ as long
as possible we still get a spine of $X$. To do so, let us
triangulate $P$ in such a way that $S(X)$ intersects the $2$-simplices only, and
each $2$-simplex at most once,
and let us examine the elementary collapses. Of course we still have a spine of
$X$ after collapsing a $1$-simplex. The only case where collapsing
a $2$-simplex $\sigma$ does not give a spine arises when $\sigma$
intersects $S(X)$, as shown in Fig.~\ref{nonnuclear:fig}-left.
    \begin{figure}
    \begin{center}
    \input{nonnuclear.pstex_t} 
    \nota{Obstructions to collapsing $2$-simplices} \label{nonnuclear:fig}
    \end{center}
    \end{figure}
Since $\sigma$ has at least one free edge, the
same component $B$ of $X\setminus P$ is incident to both sides of $\sigma$.
By definition of spine, $B$ is a discal $3$-orbifold of
ordinary or cyclic type. Depending on which of these cases occurs,
we have in $X$ either the situation shown in
Fig.~\ref{nonnuclear:fig}-centre, which implies that $X$ has an ordinary $\cyctres(p)$
connected summand, or the situation shown in
Fig.~\ref{nonnuclear:fig}-right, which implies that $X$ contains the bad $2$-suborbifold
$(S^2;q)$. A contradiction.

We have shown so far that any spine $P$ of $X$ can be replaced by the
\emph{nuclear} spine obtained by collapsing $P$ as long as possible.
The conclusion now follows from the obvious remark that $X$ does have spines
$P$ with $c(P,S(X))=c(X)$, and from the following assertions:
\begin{enumerate}
\item If $P$ is a nuclear polyhedron and $Q\subset P$ is
homeomorphic to $P$ then $Q$ is equal to $P$;
\item If $\{P_i\}$ is a sequence of nuclear polyhedra with $P_{i+1}\subseteq P_i$
then $P_i$ is eventually constant.
\end{enumerate}
To prove (1), note first that $S(Q)\subset S(P)$. It easily
follows that $S(Q)=S(P)$, otherwise $Q$ would be collapsible, but
$Q\cong P$. We now claim that the $2$-dimensional portion of $Q$
coincides with that of $P$. By contradiction, let $\alpha$ be a
path in the $2$-dimensional portion of $P$ which joins a point of
$S(P)=S(Q)$ to a point of $P\setminus Q$. If $x$ is the last point
of $\alpha$ which belongs to $Q$, then  $x$ belongs to some face
of $Q$ which can be collapsed, so $Q$ is not nuclear, but $Q\cong
P$. The proof that $Q$ has the same $1$-dimensional portion as $P$
is carried out along the same lines.

Assertion (2) is already implicit in Matveev's statement that minimal
spines for manifolds exist, so we will refrain from giving a formal proof.
\end{proof}

\paragraph{Some non-special spines}
Besides the orbifold $\matS^3_*$ already defined above, we will
need to consider in the sequel certain orbifolds $(\matP^3,F_p)$
and $(L_{3,1},F_p)$. In both cases $F_p$ is a circle of order $p$,
given by a non-singular fibre of the natural Seifert fibration. To
include the case of the manifolds $S^3,\matP^3,L_{3,1}$ (those
shown by Matveev to be the closed irreducible ones of complexity
$0$) we stipulate that if $K$ is a knot then $K_p$ denotes $K$
equipped with the cone order $p$ if $p\geqslant 2$, and it denotes
the empty set if $p=1$. Coherently with this choice, we denote
$\ordtres$ also by $\cyctres(1)$.

\begin{prop}
$c(\cyctres(p))=c(\matP^3,F_p)=c(L_{3,1},F_p)=p-1,$ and any minimal spine of
any of these orbifolds is non-special.
\end{prop}

\begin{proof}
The second assertion follows from the first one, because if $P$ is
a special spine of $(M,K_p)$, where $K$ is a knot, then $P$ has at
least one vertex and it meets $K$ at least once, so
$c(P,K_p)\geqslant 1+(p-1)>p-1$. For the same reason it is
sufficient to show that for each of the orbifolds in question
there is a spine which has no vertices and meets the singular set
once. The $2$-disc, the projective plane, and the ``triple hat''
$$\Big\{z\in\matC:\ |z|\leqslant 1\Big\}\Big/{\Big(z\sim w\ {\rm if}\ |z|=|w|=1\ {\rm and}\ z^3=w^3\Big)}$$
are such spines for
$\cyctres(p)$, $(\matP^3,F_p)$, and $(L_{3,1},F_p)$, respectively.
\end{proof}

The next result deals with the order-$(p,q,r)$ vertex spherical $3$-orbifold $\vertres(p,q,r)$
introduced in Section~\ref{def:section}.

\begin{prop}\label{theta:prop}
$c(\vertres(p,q,r))=(p-1)+(q-1)+(r-1)$ and any minimal
spine of this orbifold is not special.
\end{prop}

\begin{proof}
Denote the singular set of $\vertres(p,q,r)$ by $\theta_{p,q,r}$.
A sphere $S^2$ embedded in $S^3$ so to separate one vertex of $\theta_{p,q,r}$ from
the other one is a spine of the required complexity.  The proof
is completed along the lines of the previous one.
\end{proof}

\paragraph{Orbifolds with minimal special spines}
We prove now the main result of the present section:

\begin{teo}\label{irred:teo}
If $X$ is an irreducible $3$-orbifold and not
$\cyctres(p)$, $(\matP^3,F_p)$, $(L_{3,1},F_p)$, or $\vertres(p,q,r)$, then any minimal spine
of $X$ is special.
\end{teo}

\begin{proof}
Let $P$ be a minimal spine of $X$.
By Proposition~\ref{nuclear:prop} we can assume that $P$ is nuclear.
Then $P$ is the union of an almost-special polyhedron and a graph
(see~\cite{Matveev:AAM}), so it can be non-special only if one of the following occurs:
\begin{enumerate}
\item $P$ is a point; \item $P$ has some purely $1$-dimensional
portion; \item $P$ is the triple hat; \item $P$ is almost special
and some $2$-dimensional region of $P$ is not a $2$-disc.
\end{enumerate}
In case (1) of course $X$ is $\ordtres$. Suppose we are in case
(2), and take a small $2$-disc $D$ transversal to an isolated edge
of $P$. The circle $\partial D$ can be seen as a loop on the
boundary of some component $B$ of $X\setminus P$. Recall that $B$
is one of the $\matD^{\,3}_*$'s. The condition that $X$ contains
no bad $2$-suborbifolds easily implies that $\partial D$, viewed
on $\partial B$, cannot encircle any of the cone points of
$\partial B$. Therefore $\partial D$ bounds within $B$ a
non-singular $2$-disc $D'$. Now $D\cup D'$ is an ordinary sphere,
so it bounds an ordinary $3$-disc $\Delta$, and by construction
$\partial\Delta$ meets $P$ in a single point. It follows that if
we dismiss the whole of $P\cap \Delta$ we still have a spine of
$X$. Since $P$ is nuclear and $P\setminus\Delta$ is not, we get a
contradiction to the minimality of $P$.

In case (3) it is easy to check that $|X|$ is necessarily
$L_{3,1}$. Moreover every point of intersection between $S(X)$ and
$P$ gives two cone points on the boundary of the $3$-disc
$X\setminus P$, so there is at most one intersection point. It
easily follows that $X$ is $(L_{3,1},F_p)$ for some $p$.

Let us turn to case (4). One can see that one of the following must occur:
\begin{enumerate}
\item[(a)] $P=S^2$; \item[(b)] $P=\matP^3$; \item[(c)] there
exists a 2-dimensional region $R$ of $P$ and a loop $\gamma$ on
$R$ which cuts $R$ into two surfaces neither of which is a
$2$-disc.
\end{enumerate}
In case (a) we see that $X$ is either $\cyctres(p)$ or
$\vertres(p,q,r)$. In case (b) we have $|X|=\matP^3$ and, with the
same argument used for $L_{3,1}$, we see that $X=(\matP^3,F_p)$.
Let us then turn to case (c), and note that $\gamma$ is
orientation-preserving on $R$, so $R$ is transversely orientable
along $\gamma$ (because $|X|$ is orientable).  It follows that
there exists an annulus $A$ with core $\gamma$ which meets $P$
along $\gamma$ only, and transversely. The boundary components
$\gamma_\pm$ of $A$ can now be viewed as loops on components
$B_\pm$ of $X\setminus P$, possibly with $B_-=B_+$. We take now
within $B_\pm$ a disc $D_\pm$ bounded by $\gamma_\pm$, with
minimal possible intersection with $S(X)$, and $D_-\cap D_+ $
empty in case $B_-=B_+$. We then get a sphere $\Sigma=D_+\cup
A\cup D_-$. Recall again that on $\partial B_\pm$ there can be
zero, two, or three cone points, and note that $\gamma_\pm$ could
or not separate one of these cone points from the other one(s).
Depending on which case occurs for $\gamma_-$ and $\gamma_+$, the
sphere $\Sigma$ can be either non-singular, or a bad
2-suborbifold, or a spherical $2$-suborbifold of type
$\cycdues(p)$. By the irreducibility of $X$ we deduce that the
second case is impossible, and in the other cases there is an
ordinary or cyclic discal $3$-orbifold $\Delta$ bounded by
$\Sigma$. By construction $\partial\Delta\cap P=\gamma$ and
$P\cap\Delta$ is not a $2$-disc. If $\Delta$ is ordinary, we
conclude that $P$ is non-minimal replacing $P\cap\Delta$ by a
$2$-disc bounded by $\gamma$, as in~\cite{Matveev:AAM}, because
either $P\cap\Delta$ contains vertices or $P\cap\Delta$ properly
contains such a $2$-disc. Suppose then that $\Delta$ is cyclic,
with singular set an arc $\alpha$. By Lemma~\ref{arc:lem},
$\alpha\cap P$ consists of at most one point. If $\alpha\cap P$ is
empty then we get from $P$ a new spine $P'$ simply by dismissing
the whole of $P\cap\Delta$. Similarly, if $\alpha\cap P$ consists
of a point, we obtain from $P$ a new spine $P'$ by replacing
$P\cap\Delta$ by a $2$-disc bounded by $\gamma$. Note that in both
cases $P'$ has the same intersections as $P$ with $S(X)$. Moreover
$P'$ contradicts the minimality of $P$ for precisely the same
reasons as in the case where $\Delta$ is an ordinary $3$-disc. The
proof is complete.
\end{proof}

\begin{cor}
For every natural $n$ the set
$$\calX_n:=\big\{X\ \textrm{irreducible}\ 3{-orbifold}:\ c(X)\leqslant n\big\}$$
is finite.
\end{cor}

\paragraph{Duality}
We define a \emph{triangulation} of a $3$-orbifold $X$ to be a triangulation
of $|X|$ such that $S(X)$ is a subset of the $1$-skeleton.  Recall that for manifolds
we use a notion of triangulation more flexible than the traditional PL notion
(see Section~\ref{def:section}). Note also that we do not impose the restriction,
used for instance in~\cite{BMP}, that each tetrahedron should have at most one singular edge.
We now extend Proposition~\ref{duality:basic:prop}, showing that for
irreducible orbifolds minimal spines typically are dual to triangulations.

\begin{prop}\label{duality:orbi:prop}
Under the assumptions of Theorem~\ref{irred:teo}, dual to any minimal spine of $X$ there is
a triangulation.
\end{prop}

\begin{proof}
We know that $P$ is special, so dual to it there is a
triangulation $T$ of $|X|$, by
Proposition~\ref{duality:basic:prop}. We now claim that every
region $R$ of $P$ meets $S(X)$ at most once. This claim easily
yields the conclusion, because we can then isotope the arc of
$S(X)$ meeting $R$, if any, to the edge of $T$ dual to $R$.

To prove the claim, suppose by contradiction that some region $R$ has $k\geqslant 2$
intersections. Of course $k\leqslant 3$. In case $k=2$ we see that the two components
$B_\pm$ of $X\setminus P$ incident to $R$ must be distinct, because the total number of cone
points on their boundaries is at least $4$. By definition of spine, $S(B_\pm)$ is either
a trivial arc or a trivial Y-graph. This leads to the three possibilities
shown in Fig.~\ref{onecrossing:fig}.
    \begin{figure}
    \begin{center}
    \mettifig{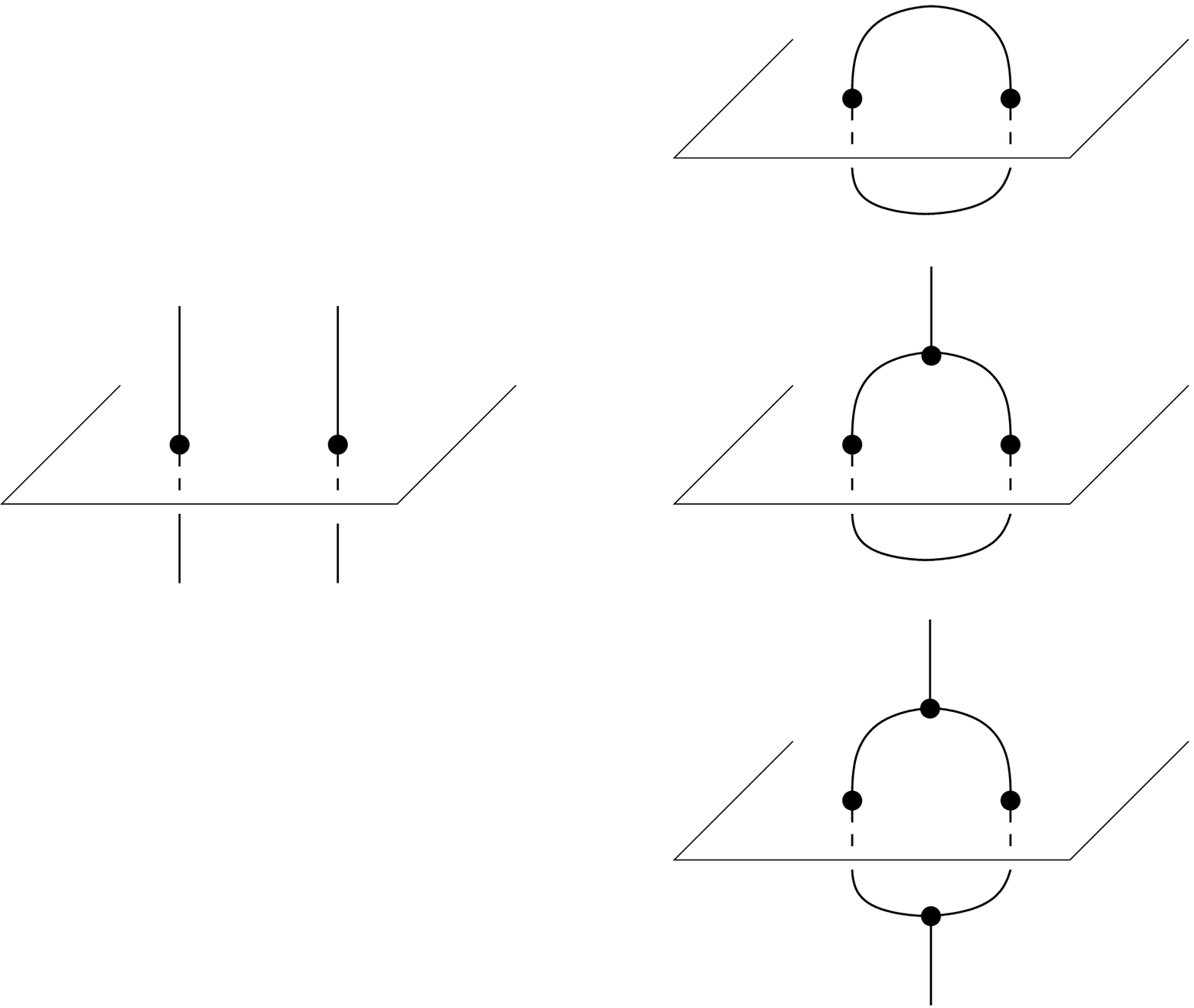,width=8cm}
    \nota{A region of spine meeting the singular locus twice.} \label{onecrossing:fig}
    \end{center}
    \end{figure}
Correspondingly, one sees either that $X$ has an $\cyctres(p)$ ordinary connected summand,
or that $X$ contains a bad $2$-suborbifold, or that $X$ has a cyclic connected summand $\vertres(p,q,r)$.
All these cases are excluded by the assumptions. Similarly,
if $k=3$ then $X$ has an ordinary connected summand $\vertres(p,q,r)$, whence the conclusion.
\end{proof}

\section{Spherical splitting of 3-orbifolds}\label{splitting:section}
To prove Theorems~\ref{additivity:intro:teo} and~\ref{estimates:intro:teo}
we need to recall more about connected sum of orbifolds than we did in
Theorem~\ref{splitting:cite:teo}. We address the reader to~\cite{orb:split}
for all proofs and more details.

\paragraph{Splitting systems}
There is a dual viewpoint to connected sum, which we will mostly employ.
To describe it, we first give a definition.
If a $3$-orbifold $Y$ is bounded by spherical $2$-orbifolds, we can canonically
associate to $Y$ a closed orbifold $\widehat{Y}$ by attaching the appropriate discal
$3$-orbifold to each component of $\partial Y$. We say that $\widehat{Y}$ is obtained
by \emph{capping} $Y$.

Recall now that an operation of connected sum consists of two
phases, namely puncturing and gluing. If we consider now a
sequence of connected sums, we can always arrange the successive
punctures to be disjoint from the spherical $2$-orbifolds along
which the previous gluings have been performed. If $X$ is the
result of the sequence of sums, we then have in $X$ a finite
system $\calS$ of spherical $2$-suborbifolds such that the
components of $(X\setminus U(\calS))\widehat{\ }$ are the original
orbifolds we have summed. (Just as above, $U$ denotes here an open
regular neighbourhood, and, for the sake of brevity, we will
actually identify $X\setminus U(\calS)$ to $X\setminus \calS$).

Conversely, to any system $\calS$ of \emph{separating} spherical
$2$-orbifolds in a $3$-orbifold $X$ there corresponds a
realization of $X$ as a connected sum of the components of
$(X\setminus \calS)\widehat{\ }$. We now call \emph{efficient} a
finite system of separating spherical $2$-suborbifolds which
corresponds to an efficient connected sum. To characterize the
efficient splitting systems, we call \emph{punctured discal} a
3-orbifold obtained from one of the $\matD^3_*$'s by removing a
regular neighbourhood of a finite subset.

\begin{prop}
A system $\calS$ of spherical $2$-suborbifolds of a $3$-orbifold $X$ is efficient if and only if no component
of $X\setminus \calS$ is punctured discal, and each component of $(X\setminus \calS)\widehat{\ }$
is irreducible.
\end{prop}

\paragraph{Uniqueness and stepwise splitting}
The core of the uniqueness part of Theorem~\ref{splitting:cite:teo}
is the following result:

\begin{prop}
Let $\calS$ be an efficient splitting system in a closed $3$-orbifold $X$.
Then $\calS$ and $(X\setminus \calS)\widehat{\ }$, as abstract collections
of $2$- and $3$-orbifolds respectively, depend on $X$ only.
\end{prop}

We now state a result which gives a practical recipe to construct an efficient
splitting system. Recall that a spherical $2$-orbifold is inessential if it bounds
a discal $3$-orbifold.

\begin{teo}\label{step:split:teo}
Let $X$ be a closed locally orientable $3$-orbifold. Suppose that $X$ does not contain
any bad $2$-suborbifold, and that every spherical $2$-suborbifold of $X$ is separating.
Starting with $\calS=\emptyset$, consider the (non-deterministic) process
described by the following steps:
\begin{itemize}
\item[1.] If all the ordinary spherical $2$-suborbifolds of $(X\setminus \calS)\widehat{\ }$
are inessential, turn to Step 2. Otherwise choose $\Sigma$ as one such $2$-suborbifold disjoint from $\calS$,
redefine $\calS$ as the given $\calS$ union $\Sigma$, and repeat Step 1;
\item[2.] If all the cyclic spherical $2$-suborbifolds of $(X\setminus \calS)\widehat{\ }$ are inessential, turn to Step 3.
Otherwise choose $\Sigma$ as one such $2$-suborbifold disjoint from $\calS$, redefine $\calS$
as the given $\calS$ union $\Sigma$, and repeat Step 2;
\item[3.] If all the vertex spherical $2$-suborbifolds of $(X\setminus \calS)\widehat{\ }$ are inessential, output $\calS$.
Otherwise choose $\Sigma$ as one such $2$-suborbifold disjoint from $\calS$,
redefine $\calS$ as the given $\calS$ union $\Sigma$, and repeat Step 3.
\end{itemize}
Then the process is finite and the final $\calS$ is an efficient splitting system.
\end{teo}

\section{Handles and normal 2-orbifolds}\label{handles:section}
In this section we introduce the notion of handle decomposition
for a $3$-orbifold, we define normal $2$-orbifolds with respect to
handle decompositions, and we prove (under suitable assumptions)
the fundamental fact that essential spherical 2-orbifolds can be
normalized.

\paragraph{Handle decompositions of 3-orbifolds}
We will use Matveev's terminology of~\cite{Matveev:ATC3M}, calling
\emph{balls}, \emph{beams}, and \emph{plates} respectively the
$0$-, $1$-, and $2$-handles of a handle decomposition of a
$3$-manifold $M$. In addition, we will call \emph{caps} the
3-handles. Again following Matveev, we will also call
\emph{islands} (respectively, \emph{bridges}) the connected
components of the attaching loci between the balls and the beams
(respectively, the balls and the plates).

If $X$ is a 3-orbifold, we call
\emph{orbifold-handle decomposition} of $X$ a handle decomposition of $|X|$ such
that the following holds:
\begin{itemize}
\item The balls and the beams are disjoint from $S(X)$;
\item Each plate $D^2\times I$ meets $S(X)$ in an arc $\{*\}\times I$;
\item Each cap is isomorphic to some $\matD^{\,3}_*$.
\end{itemize}

Existence of orbifold-handle decompositions is very easily established.

\paragraph{Normal 2-orbifolds}
Fix an orbifold-handle
decomposition of a $3$-orbifold $X$. We will say that a 2-suborbifold
$\Sigma$ of $X$ is in \emph{normal position} with respect to the decomposition if:
\begin{itemize}
\item $\Sigma$ is disjoint from the caps; \item $\Sigma$ meets
each plate $D^2\times I$ in a family of parallel $2$-discs
$D^2\times\{t_1,\ldots,t_\nu\}$; \item $\Sigma$ meets each beam
$I\times D^2$ in a set of the form
$I\times(\alpha_1\cup\ldots\cup\alpha_\mu)$, where $\{\alpha_j\}$
is a family of disjoint properly embedded arcs in $D^2$; \item
$\Sigma$ meets each ball in a union of $2$-discs, and the boundary
loop of each of these $2$-discs passes through each bridge at most
once.
\end{itemize}

We prove now the key result of the section.
In the proof we will refer to the ``normalization'' moves for surfaces
described in~\cite{Matveev:ATC3M}.

\begin{prop}\label{normalization:prop}
Let $X$ be a closed $3$-orbifold which contains no bad $2$-suborbifolds
and no non-separating spherical suborbifold. Fix an
orbifold-handle decomposition of $X$. Then:
\begin{enumerate}
\item If $X$ contains some essential ordinary sphere then it contains a normal one;
\item Suppose that in $X$ every ordinary sphere is inessential. If $X$ contains an essential
cyclic or vertex spherical $2$-suborbifold then it contains a normal one of
the same type.
\end{enumerate}
\end{prop}

\begin{proof}
We prove both points at the same time.
Let $\Sigma$ be the given essential
spherical $2$-suborbifold. We note that $|\Sigma|$ is a sphere and we
apply to it the normalization moves
for surfaces, noting that all of them except two
are isotopies of $\Sigma$ (as a 2-suborbifold of $X$, not just as a surface), so they preserve
essentiality.  The two exceptions are as follows:
\begin{itemize}
\item The compression of $|\Sigma|$ along a $2$-disc, followed by
the choice of one of the two resulting spheres; \item The move
which allows to eliminate double passages along bridges, as shown
in Fig.~\ref{doublebridge:fig} (and explained in its caption).
\end{itemize}
    \begin{figure}
    \begin{center}
    \mettifig{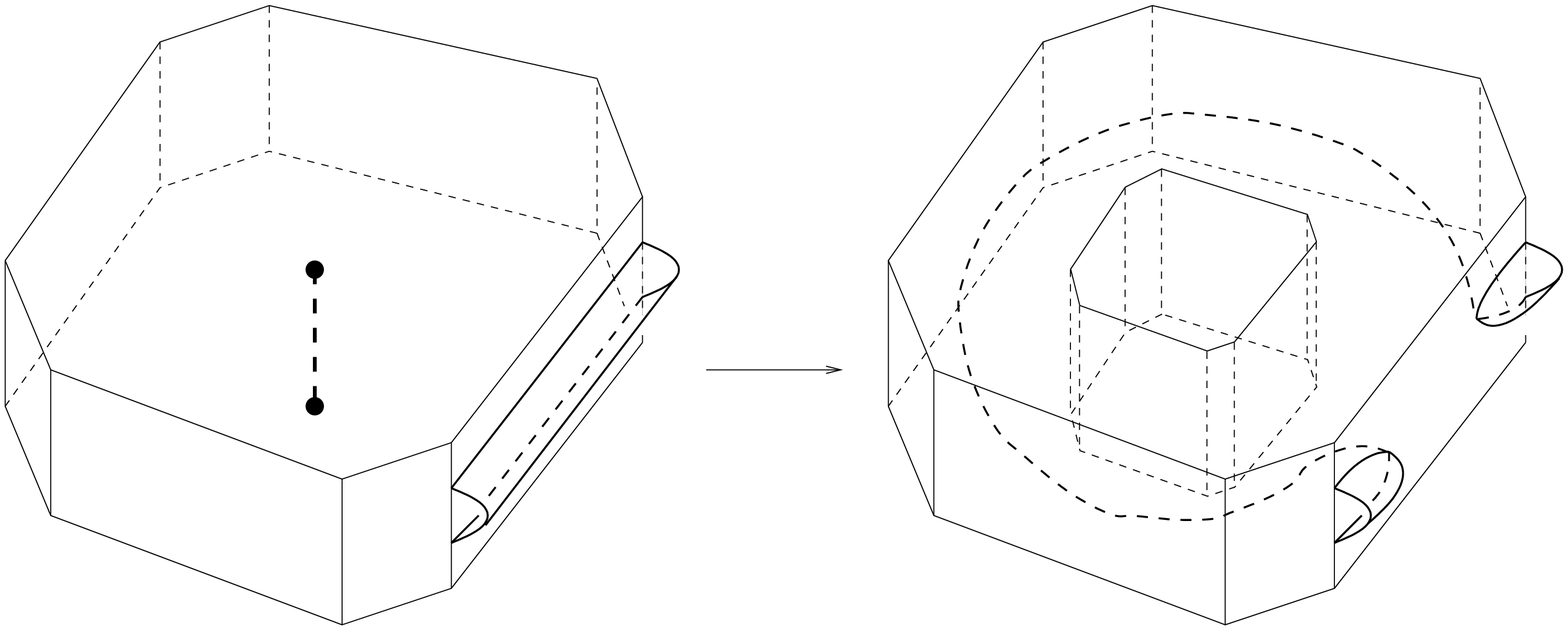,width=10cm}
    \nota{To eliminate a double passage of a surface along a bridge one should first isotope the surface
    so that it avoids the core of the corresponding plate, and then
    resize the plate to a neighbourhood of its core.} \label{doublebridge:fig}
    \end{center}
    \end{figure}
Since the first move takes place within the union of balls and
beams, the compression $2$-disc involved does not intersect the
singular set. An equivalent way of describing the second move is
provided in Fig.~\ref{alsocompr:fig}:
    \begin{figure}
    \begin{center}
    \input{alsocompr.pstex_t} 
    \nota{Assuming the core of the plate has order $p\geqslant1$, the move of Fig.~\ref{doublebridge:fig}
    can be seen as the compression shown here followed by the elimination of the
    inessential ordinary or cyclic spherical $2$-orbifold $\Sigma'$.} \label{alsocompr:fig}
    \end{center}
    \end{figure}
we first compress $\Sigma$ along the $2$-disc $D$ and then we
dismiss the resulting $2$-orbifold $\Sigma'$, which bounds an
ordinary or cyclic discal $3$-orbifold, keeping the $2$-orbifold
$\Sigma''$.

According to this discussion, the following claim is sufficient to
conclude the proof: \emph{If $\Sigma'$ and $\Sigma''$ are the
$2$-suborbifolds obtained by compressing $\Sigma$ along a $2$-disc
$D$ which does not meet $S(X)$, then either $\Sigma'$ or
$\Sigma''$ is essential and has the same nature as $\Sigma$.}

Let $D'$ and $D''$ be the closures of the components of $\Sigma\setminus\partial D$. Since $D$ is non-singular, the
assumption that $X$ contains no bad $2$-suborbifold easily implies that (up to changing notation)
also $D''$ is non-singular. Therefore $\Sigma'=D\cup D'$ has the same nature as $\Sigma$
and $\Sigma''$ is ordinary, and the conclusion readily follows if we show that either
$\Sigma'$ or $\Sigma''$ is essential.

Suppose by contradiction that both $\Sigma'$ and $\Sigma''$ are inessential,
\emph{i.e.} that they bound discal $3$-orbifolds.  Then $\Sigma$ bounds in $X$
an orbifold of one of the following two types:
\begin{itemize}
\item The gluing along a non-singular $2$-disc of two discal $3$-orbifolds, one of which is ordinary;
\item One of the components of a $3$-orbifold obtained by cutting a discal $3$-orbifold along
a non-singular $2$-disc.
\end{itemize}
Any orbifold of such a type is discal, and the proof is complete.
\end{proof}

\paragraph{Normal 2-orbifolds with respect to spines}
Suppose $P$ is a simple spine of a $3$-orbifold
$X$, and consider a triangulation of $P$,
in the traditional PL sense, with the property that $1$-simplices
are disjoint from $S(X)$ and each 2-simplex meets $S(X)$ in at most one point.
Note that the last condition may force subdividing some $2$-dimensional portion
of $P$ more than a triangulation of $P$ itself would require, because the $2$-dimensional
components of $P$ can a priori intersect $S(X)$ in as many as three points.
Consider the handle decomposition of $|X|$ whose balls, beams, and plates are obtained by
thickening the vertices, edges, and triangles of this
triangulation of $P$, and taking as caps the components of what is left.
One easily sees that this is also an orbifold-handle decomposition of $X$.
Now let $\Sigma$ be a normal $2$-orbifold with respect to this orbifold-handle
decomposition. Then we can attach to each 2-dimensional component
of $P$ an integer, corresponding to the number of times $\Sigma$ runs parallel
to the component, and $\Sigma$ is uniquely determined by the assignment of
non-negative integers to the 2-components of $P$.

\section{Additivity under ordinary\\ and cyclic connected sum}\label{additivity:section}
In this section we deal with two of the three types of connected
sum, leaving the vertex type to the next section.  We show that
(under some restrictions) orbifold complexity is additive (in a
suitable sense) under these operations.

\paragraph{Subadditivity}
We begin by giving upper estimates for the complexity of an orbifold obtained
as an ordinary or cyclic connected sum of two given orbifolds. For cyclic connected sums
it turns out that
the estimate depends also on whether the singular components involved in the sum
both contain vertices or not. And we will see below that this distinction is
unavoidable.

\begin{prop}\label{ord:sum:prop}
If $X$ is the ordinary connected sum of $X_0$ and $X_1$ then
$$c(X)\leqslant c(X_0)+c(X_1).$$
\end{prop}

\begin{proof}
Let $P_j$ be a minimal spine of $X_j$. To perform the connected
sum, we must remove from $X_j$ an ordinary $3$-disc $\Delta_j$
disjoint from $S(X_j)$, choose a homeomorphism $h$ between
$\partial \Delta_0$ and $\partial \Delta_1$, and glue
$X_0\setminus \Delta_0$ to $X_1\setminus \Delta_1$ along $h$.
Since the position of $\Delta_j$ is immaterial, we choose
$\Delta_j$ to be disjoint from $P_j$. Therefore $P_0\sqcup P_1$
can be viewed as a subset of $X$, and its complement consists of
some admissible $3$-discs together with one component
$S^2\times(-1,1)$ which contains at most two portions of singular
set, separated by $S^2\times\{0\}$ if there are two of them. The
idea is now to construct a spine of $X$ by adding $\big(\{{\rm
*}\}\times(-1,1)\big)\cup\big( S^2\times\{0\}\big)$ to $P_0\sqcup
P_1$, see Fig.~\ref{spineofsum:fig}
    \begin{figure}
    \begin{center}
    \input{spineofsum.pstex_t} 
    \nota{Spine of an ordinary connected sum.} \label{spineofsum:fig}
    \end{center}
    \end{figure}

To formalize this idea, we observe that $\Delta_j$ is contained in
a component of $X_j\setminus P_j$, so we can find within this
component an arc $\gamma_j$ which joins $P_j$ to $\partial
\Delta_j$ without meeting $S(X_j)$. We define now
$P'_j=P_j\cup\gamma_j\cup\partial \Delta_j$ and we construct $P$
as the gluing of $P'_0$ and $P'_1$ under the homeomorphism $h$
between $\partial \Delta_0$ and $\partial \Delta_1$. It is easy to
see that $P$ is a spine of $X$ and
$c(P,S(X))=c(P_0,S(X_0))+c(P_1,S(X_1))$, whence the conclusion.
\end{proof}

\begin{rem}\emph{As opposed to the case of manifolds treated in~\cite{Matveev:AAM},
in the previous proof it was essential to include in the spine $P$
of $X$ the sphere along which the connected sum was performed.
This is because the component of $X_j\setminus P_j$ containing
$\Delta_j$ can be a cyclic or vertex $3$-disc. And, if this
happens for $j=0$ and $j=1$, the gluing of these components along
a small boundary $2$-disc does not have the required type.}
\end{rem}

\begin{prop}\label{cyc:cyc:sum:prop}
Let $X$ be a cyclic connected sum of $X_0$ and $X_1$ along arcs of order $p$.
Suppose that at least one of the singular components involved in the sum is a knot.
Then
$$c(X)\leqslant c(X_0)+c(X_1)-(p-1).$$
\end{prop}

\begin{proof}
Let $P_j$ be a minimal spine of $X_j$. At the level of the
singular set, the connected sum is performed by removing a small
arc $\alpha_0$ from $S(X_0)$ and one $\alpha_1$ from $S(X_1)$, and
by joining together in pairs the ends thus created. Of course we
can assume that $\alpha_j$ is disjoint from $P_j$, so it is
contained in a component $B_j$ of $X_j\setminus P_j$. Let
$\widetilde\alpha_j$ be the component of $S(X_j)\setminus\{{\rm
vertices}\}$ which contains $\alpha_j$. By assumption, up to
switching indices, $\widetilde\alpha_0$ is a knot. Thus, by
Lemma~\ref{arc:lem}, $\widetilde\alpha_0$ is contained in the
closure of $B_0$, and it intersects $P_0$ in a single point $x_0$.
The same lemma implies that $\widetilde\alpha_1$ also meets $P_1$
in a single point $x_1$. A schematic representation of the
situation is given in Fig.~\ref{subadd:fig}-left, assuming that
$B_1$ is a vertex $3$-disc.
    \begin{figure}
    \begin{center}
    \input{subadd.pstex_t} 
    \nota{Construction of a spine of a cyclic connected sum.} \label{subadd:fig}
    \end{center}
    \end{figure}
Note that the component $B_2$ of $X_1\setminus P_1$ may or not be equal to $B_1$.

Without loss of generality we can now suppose that the connected
sum is performed by removing $3$-discs contained in $B_0$ and
$B_1$. Then $P_0\sqcup P_1$ can be viewed as a subset of $X$, and
its complement consists of some admissible $3$-discs together with
a component $S^2\times(-1,1)$.  This component contains two
portions of singular set, namely an arc which joins
$S^2\times\{-1\}$ to $S^2\times\{+1\}$, and another portion which
can be an arc or a Y-graph and has precisely one end on
$S^2\times\{-1\}$.

We can now create a spine $P$ of $X$ by adding to $P_0\sqcup P_1$ a cylinder
in $S^2\times(-1,1)$ which encircles the singular arc, as suggested in
Fig.~\ref{subadd:fig}-right. Note that $c(P,S(X))=c(P_0,S(X_0))+c(P_1,S(X_1))$.
Let us consider now the components of $X\setminus P$, denoted
by $\widetilde B_j$ in the picture. Then of course $\widetilde B_3$ is not equal
to $\widetilde B_1$ or $\widetilde B_2$ (but we could have $\widetilde B_1=\widetilde B_2$
in case $B_1=B_2$). Therefore we can puncture $P$ near $x_0$ or near $x_1$ (but not both, in
general), thus reducing by $1$ the number of intersection points between $P$ and $S(X)$.
The conclusion readily follows.
\end{proof}

\begin{prop}\label{cyc:vert:sum:prop}
Let $X$ be a cyclic connected sum of $X_0$ and $X_1$.
Suppose that both the singular components involved in the sum contain vertices.
Then
$$c(X)\leqslant c(X_0)+c(X_1).$$
\end{prop}

\begin{proof}
In the previous proof, where did we use the assumption that the
connected sum should involve at least one singular component
without vertices? It was when we supposed that the $3$-disc $B_0$
which contains $\alpha_0$ intersects $S(X_0)$ in an arc (rather
than in a Y-graph). This condition need not be met here, but we
can always apply to $P_0$ a bubble move near a point of
intersection between $P_0$ and $S(X_0)$. This creates a new spine
$P'_0$ of $X_0$ with $c(P'_0,S(X_0))= c(P_0,S(X_0))+(p-1)$. We can
now apply to $P'_0$ and $P_1$ precisely the same argument as
above, deducing that $c(X)\leqslant
c(P'_0,S(X_0))+c(P_1,S(X))-(p-1)$, whence the conclusion.
\end{proof}

\paragraph{Superadditivity}
In this paragraph we prove the estimates opposite to those established in
Propositions~\ref{ord:sum:prop},~\ref{cyc:cyc:sum:prop}, and~\ref{cyc:vert:sum:prop}.

\begin{prop}\label{inv:ord:sum:prop}
Let $X$ be a $3$-orbifold. Let
$\Sigma$ be a separating ordinary spherical $2$-suborbifold of $X$. Let
$X_+$ and $X_-$ be the components of $(X\setminus\Sigma)\widehat{\ }$.
Let $P$ be a spine of $X$.
Suppose that $\Sigma$ is in normal position with respect to $P$.
Then there exist spines $P_\pm$ of $X_\pm$ such that
$$c(P_+,S(X_+))+c(P_-,S(X_-))\leqslant c(P,S(X)).$$
\end{prop}

\begin{proof}
The desired spines are obtained by
cutting $P$ along $\Sigma$, as in~\cite{Matveev:AAM}.
As we cut we do not create vertices, so $c(P_+)+c(P_-)\leqslant c(P)$.
Recall now that a region of $P$ having a positive weight $n$ in $\Sigma$
gets replaced by $n$ regions, so the intersections between the spine and the singular
set could a priori increase after cutting. However, our sphere $\Sigma$ was supposed to
be ordinary, so the regions carrying positive weight do not intersect $S(X)$,
and the conclusion easily follows.
\end{proof}

\begin{prop}\label{inv:cyc:cyc:sum:prop}
Let $X$ be a $3$-orbifold. Let
$\Sigma$ be an order-$p$ cyclic separating spherical $2$-suborbifold of $X$.
Suppose there exists an arc in $S(X)$ joining the two cone points of $\Sigma$ and
not containing vertices of $S(X)$.
Let $X_+$ and $X_-$ be the components of $(X\setminus\Sigma)\widehat{\ }$.
Let $P$ be a minimal spine of $X$.
Suppose that $\Sigma$ is in normal position with respect to $P$.
Then there exist spines $P_\pm$ of $X_\pm$ such that
$$c(P_+,S(X_+))+c(P_-,S(X_-))\leqslant c(P,S(X))+(p-1).$$
\end{prop}

\begin{proof}
As in the previous proof, we cut $P$ along $|\Sigma|$, but now we must be careful because
the singular set indeed intersects the spine.  To be precise, we note that $\Sigma$ meets $S(X)$ twice, so we
have the following possibilities:
\begin{itemize}
\item[(a)] There is a region of $P$ which intersects $S(X)$ once and has weight $2$ in $\Sigma$;
all other regions of $P$ having positive weight do not intersect $S(X)$;
\item[(b)] There is a region of $P$ which intersects $S(X)$ twice and has weight $1$ in $\Sigma$;
all other regions of $P$ having positive weight do not intersect $S(X)$;
\item[(c)] There are two regions regions of $P$ which both intersect $S(X)$ once and have weight $1$ in $\Sigma$;
all other regions of $P$ having positive weight do not intersect $S(X)$.
\end{itemize}
However we can see that cases (b) and (c) are absurd, because, together with the assumption
that the cone points of $\Sigma$ can be joined in $S(X)$ avoiding vertices,
they would imply that there is an arc of $S(X)$ meeting $P$ twice and not containing vertices,
which contradicts Lemma~\ref{arc:lem} (note that we have supposed $P$ to be minimal).

The situation is therefore as shown (in a cross-section) in Fig.~\ref{cycut:fig}-left, where
    \begin{figure}
    \begin{center}
    \input{cycut.pstex_t} 
    \nota{Cutting a spine along a sphere which passes twice through a region.} \label{cycut:fig}
    \end{center}
    \end{figure}
$\Sigma'$ and $\Sigma''$ denote two of the spherical $2$-orbifolds bounding a regular neighbourhood of $P$
(possibly $\Sigma'=\Sigma''$). After cutting we have then the situation of Fig.~\ref{cycut:fig}-right, where the
new spines are denoted by $\widetilde P_\pm$, and $\Sigma_\pm$ are the spherical orbifolds
bounding a product neighbourhood of $\Sigma$. (To draw the figure we have used the assumption that $\Sigma$ is
separating). We note now that in $X_+$ we certainly have $\Sigma_+\neq \Sigma'$, so we can puncture at least one region of
$\widetilde P_+$, eliminating one intersection point with $S(X_+)$. We define $P_+$ as the punctured $\widetilde P_+$,
and $P_-$ as $\widetilde P_-$.  We note that $P_+$ is still a spine of $X_+$ because (by construction)
$\Sigma_+$ bounds a cyclic discal $3$-orbifold in $X_+$.
Since in passing from $P$ to $P_+\sqcup P_-$ we have not created vertices
but we have doubled a point of intersection with the singular locus, the desired estimate readily follows.
\end{proof}

\begin{prop}\label{inv:cyc:vert:sum:prop}
Let $X$ be a $3$-orbifold. Let
$\Sigma$ be an order-$p$ cyclic separating spherical $2$-suborbifold of $X$.
Suppose that any arc in $S(X)$ joining the two cone points of $\Sigma$
contains some vertex of $S(X)$.
Let $X_+$ and $X_-$ be the components of $(X\setminus\Sigma)\widehat{\ }$.
Let $P$ be a spine of $X$.
Suppose that $\Sigma$ is in normal position with respect to $P$.
Then there exist spines $P_\pm$ of $X_\pm$ such that
$$c(P_+,S(X_+))+c(P_-,S(X_-))\leqslant c(P,S(X)).$$
\end{prop}

\begin{proof}
The proof closely imitates the previous one, except that it is case (a) that is now absurd.
We then have the situation of Fig.~\ref{vertcut:fig}-left, where the four spherical
    \begin{figure}
    \begin{center}
    \input{vertcut.pstex_t} 
    \nota{Cutting a spine along a sphere which passes once through two region.} \label{vertcut:fig}
    \end{center}
    \end{figure}
orbifolds $\Sigma_j^{(i)}$ could be anything from equal to each
other to distinct from each other. After cutting we then have the
situation of Fig.~\ref{vertcut:fig}-right, to draw which we have
again used the fact that $\Sigma$ separates. Now $\Sigma^{(i)}_1$
is distinct from $\Sigma_\pm$, and $\Sigma_\pm$ bounds a cyclic
discal $3$-orbifold in $X_\pm$, so we can puncture $\widetilde
P_\pm$ at least once, getting a new spine $P_\pm$. The desired
estimate now follows from the fact that we have created two new
intersections with the singular set while cutting $P$, but we have
then eliminated two such intersections by puncturing.
\end{proof}

\paragraph{Summarizing statement}
It is now easy to establish the main result on additivity stated in the Introduction.

\dimo{additivity:intro:teo}
By Theorem~\ref{step:split:teo}, the unique splitting of $X$ given by
Theorem~\ref{splitting:cite:teo} is obtained by first splitting
$X$ along ordinary spheres as long as possible, and then along cyclic spheres.
Proposition~\ref{normalization:prop} shows that at each step of this successive
splitting we can take the sphere to be normal with respect to some handle
decomposition. The conclusion then follows from
Propositions~\ref{ord:sum:prop} to~\ref{inv:cyc:vert:sum:prop}.
\finedimo

\section{Estimates on vertex connected sum}\label{estimates:section}
When the efficient splitting of an orbifold $X$ into irreducible ones $X_1,\ldots,X_n$ involves cutting
along spherical $2$-orbifolds of vertex type we cannot prove results as precise as Theorem~\ref{additivity:intro:teo}.
We can however provide upper and lower estimates on the complexity of $X$ in terms of that
of $X_1,\ldots,X_n$.

\paragraph{Upper estimates}
The following result gives an upper bound on the complexity after
a vertex connected sum, under assumptions general enough to
eventually prove Theorem~\ref{estimates:intro:teo}.

\begin{prop}\label{vert:sum:prop}
Let $X_0$ and $X_1$ be orbifolds, and let $X$ be a vertex connected sum of $X_0$ and $X_1$.
Let $P_j$ be a spine of $X_j$ which is either special or homeomorphic to $S^2$.
Then there exists a special spine $P$ of $X$ such that
$$c(P,S(X))\leqslant 6\cdot (c(P_0,S(X_0))+c(P_1,S(X_1))).$$
\end{prop}

\begin{proof}
Let $B_j$ be the component of $X_j\setminus P_j$ which contains the vertex of $S(X_j)$ along
which the connected sum is performed. Let $p^{(j)}_i$, for $i=1,2,3$, be the points of intersection
between $\partial B_j$ (which is contained in $P_j$) and $S(X_j)$. Note that two of these points
could coincide, but not all three of them, so we suppose $p^{(j)}_1\neq p^{(j)}_2$.
A spine $P$ of $X$ is then constructed as follows (see Fig.~\ref{vertsum:fig}):
    \begin{figure}
    \begin{center}
    \input{vertsum.pstex_t} 
    \nota{Spine of a vertex connected sum. The region labelled $n_j$ represents $n_j$
    parallel regions} \label{vertsum:fig}
    \end{center}
    \end{figure}
\begin{itemize}
\item For $j=0,1$, remove from $P_j$ two small $2$-discs
$D^{(j)}_0$ and $D^{(j)}_1$ around $p^{(j)}_1$ and $p^{(j)}_2$;
\item For $i=1,2$, add a tube running from $\partial D^{(0)}_i$ to
$\partial D^{(1)}_i$, and a meridinal $2$-disc for this tube;
\item Add a ``rectangular'' region with one edge on each tube and,
for $j=0,1$, one edge on $\partial B_j$ running from $\partial
D^{(j)}_1$ to $\partial D^{(j)}_2$.
\end{itemize}
The fact that this spine is special is easy to prove and hence left to the reader.
Let us discuss how $c(P,S(X))$ relates to $c(P_j,S(X_j))$. Let us first ignore the vertices
created on $\partial B_0$ and $\partial B_1$ when attaching the rectangle. Then we see that
in $P$ there are two more vertices than in $P_0\sqcup P_1$, and two intersections less with
the singular set. So the complexity is not increased at this stage.
We are left to give an estimate
on how many vertices are created by attaching a region to $\partial B_j$ along a simple path $\alpha_j$
which goes from $p^{(j)}_1$ to $p^{(j)}_2$. Note that vertices arise in two forms:
\begin{itemize}
\item[(A)] At the intersections between $\alpha_j$ and $S(P_j)$, as shown in Fig.~\ref{vertsum:fig};
\item[(B)] At the self-intersections of $\alpha_j$ which can occur within the regions of $P_j$
to which $B_j$ is doubly incident.
\end{itemize}
To estimate the number of vertices of type (A), we consider a
slightly smaller $3$-disc $B'_j$ inside $B_j$, and the trivalent
graph $\Gamma_j\subset\partial B'_j$ which projects to $S(X_j)$.
Note that this graph has at most $6\cdot c(P_j)$ edges. Now we can
view $\alpha_j$ inside $\partial B'_j$ and arrange it not to touch
the edges of $\Gamma_j$ twice, and not to touch consecutive edges.
Therefore the number of vertices of type (A) arising is at most
one third the number of edges of $\Gamma_j$, hence at most $2\cdot
c(P_j)$.

Turning to type (B), we observe that at most one vertex arises in each region of $P_j$.
Note now that a region touches at least two edges, that at most three regions can be incident
to each edge, and that there are $2\cdot c(P_j)$ edges. Therefore the number of regions is
at most $3\cdot c(P_j)$, and the desired estimate immediately follows.
\end{proof}

Of course in the previous proof the upper estimate could be improved a
bit, in particular by separating in $c(P_j,S(X_j))$ the contribution given by the vertices of $P_j$
from that given by the intersections with $S(X_j)$. Since we are only interested in the qualitative
fact there exists an estimate which is linear, we refrain from doing this.

\paragraph{Lower estimates} We now turn to a lower bound, which as
usual employs normal $2$-orbifolds.

\begin{prop}\label{gen:cut:prop}
Let $X$ be a $3$-orbifold as in Theorem~\ref{splitting:cite:teo}. Let
$\Sigma$ be a separating vertex spherical $2$-suborbifold $X$. Let
$X_+$ and $X_-$ be the components of $(X\setminus\Sigma)\widehat{\ }$.
Let $P$ be a spine of $X$.
Suppose that $\Sigma$ is in normal position with respect to $P$.
Then there exist spines $P_\pm$ of $X_\pm$ such that
$$c(P_+,S(X_+))+c(P_-,S(X_-))\leqslant 4 c(P,S(X)).$$
\end{prop}

\begin{proof}
Let $\widetilde P_\pm$ be obtained by cutting $P$ along $\Sigma$,
as above. Of course $\widetilde P_\pm$ is a spine of $X_\pm$, and
the operation of cutting along $\Sigma$ obviously does not
increase the total number of vertices. This operation increases by
$3$ the number of intersections with the singular set, whence the
conclusion at once.
\end{proof}

As for the upper estimates, we have not tried to give in the previous proposition an optimal lower one.

\paragraph{Summarizing result}
We can now prove the estimation result stated in the introduction.

\dimo{estimates:intro:teo} We begin by the upper estimate. Up to
reordering we can assume that the connected sums between
$X_1,\ldots,X_k$ are of vertex type, and the other ones are not.
Then $X_1,\ldots,X_k$ are different from $\cyctres$,
$(\matP^3,F_p)$, and $(L_{3,1},F_p)$, so by
Theorem~\ref{irred:teo} they have special polyhedra or spheres as
minimal spines. Then Proposition~\ref{vert:sum:prop} implies that
$$c(X_1\#\ldots\# X_k)\leqslant 6^{k-1}\big(
c(X_1)+\ldots c(X_k)\big)$$ and the conclusion easily follows from
Propositions~\ref{ord:sum:prop},~\ref{cyc:cyc:sum:prop},
and~\ref{cyc:vert:sum:prop}.

The lower estimate immediately follows from
Theorems~\ref{splitting:cite:teo} and~\ref{step:split:teo}, and
Propositions~\ref{normalization:prop},~\ref{inv:ord:sum:prop},~\ref{inv:cyc:cyc:sum:prop},~\ref{inv:cyc:vert:sum:prop},
and~\ref{gen:cut:prop}.\finedimo

\vspace{1.5 cm}

\noindent
Dipartimento di Matematica Applicata, Universit\`a di Pisa\\
Via Bonanno Pisano 25B, 56126 Pisa, Italy\\
petronio@dm.unipi.it

\end{document}